\documentclass{arSIH}

\newtheorem{theorem}{Theorem}
\newtheorem{proposition}[theorem]{Proposition}
\newtheorem{lemma}[theorem]{Lemma}
\newtheorem{corollary}[theorem]{Corollary}

\theoremstyle{definition}

\newtheorem{definitions}[theorem]{\font\=cmssi10\Definitions\bf}

\newtheorem{remark}[theorem]{\font\=cmssi10\Remark\bf}
\newtheorem{conventions}[theorem]{\font\=cmssi10\Conventions\bf}

\def\Cal{\mathcal}
\def\Eps{\hbox{\font\=cmmi10 scaled\magstep1\\char'017}\kern0.15mm}
\def\Iota{\kern.15mm\hbox{\font\=cmmi10 scaled\magstep1\\char'023}\kern0.2mm}
\def\Nu{\hbox{\font\=cmmi10 scaled\magstep1\\char'027}\kern0.25mm}
\def\uvarPi{\kern.15mm\underline{\kern-.15mm\varPi\kern-.85mm}\kern.85mm}
\def\uOmega{\kern.3mm\underline{\kern-.3mm\Omega\kern-.3mm}\kern.3mm}

\def\tauu{\char'034}
\def\fRe{\hbox{\font\=cmr9\f\kern.1mm}\roman{Re}\kern.75mm}
\def\fIm{\hbox{\font\=cmr9\f\kern.1mm}\roman{Im}\kern.65mm}
\def\vecc#1{\kern-.5mm\vec{\kern.5mm#1}}
\def\TVS{\roman{TVS}\kern0.37mm}%
\def\LCS{\roman{LCS}\kern0.37mm}%
\def\BaS{\roman{BaS}\kern0.37mm}%
\def\dimHa{{\rm dim_{_{\kern.2mm Ha}}}}
\def\rajou{{}^{}{\Cal B}_{s\,}}
\def\Linb{\Cal L\lower.7mm\hbox{\kern.1mm\font\=cmmi6\b}}
\def\LL^#1{L\kern0.15mm\raise.4mm\hbox{$^{#1}$}\kern0.15mm}

\def\dualbeta{^{\kern0.4mm\prime}_{\kern-.2mm\raise.95mm\hbox{$_{_\beta}$}}} 
\def\Nbh{\Cal N_{\font\=cmmi6\lower.15mm\hbox{\kern.1mm\bh\kern.15mm}}}
\def\Topma{\roman{{Top_{}}_{\hbox{\font\=cmr6\ma}}\kern.15mm}}
\def\co{\hbox{\font\=cmmi12\c}\kern.15mm\lower.15mm\hbox{$_{\rm o}$}}
\def\prodc{\prod{_{_{\kern-.3mm\bold c\kern.15mm}}}}
\def\vsprod_#1_#2{\prod\kern-0.3mm{}_{_{\roman{#1}\sp{#2}\,}}} 
\def\vscoprod_#1_#2{\coprod\kern-0.3mm{}_{_{\roman{#1}\sp{#2}\,}}} 
\def\expnota^#1]_#2{\,^{#1\,]{_{}}_{\roman{#2}}}} 

\def\bold#1{{\bf#1}}
\def\roman#1{{\rm#1}}
\def\limu_#1{\lim\kern-5.5mm\lower1.5mm\hbox{$_{#1}\ $}}
\def\oseoy{\raise1.9mm\hbox{\kern.5mm\font\=cmr5\o}\kern-1.7mm y}
\def\O{{}^{}\Cal O}
\def\Univ{\hbox{\font\=cmssbx10\U}} 
\def\Pows{\Cal P\kern-.4mm_s\kern.3mm}
\def\lei{      {}_{ {}^{\,\downarrow\text{\hskip-2.1mm}       }  }  \cap       }
\def\lei{\hbox{\kern.45mm$_{^\downarrow}\kern-1.280mm\cap\kern.85mm$}}

\def\inve{\lower.85mm\hbox{$^{^-}$}\kern-.5mm{}^\iota}

\def\fvalue{\hbox{\kern.2mm\font\=cmr10\\char'022\kern-.2mm}} 
\def\ffvalue{\hbox{\kern.2mm\font\=cmr7\\char'022\kern-.2mm}} 
\def\image{\hbox{\font\=cmr10\\char'022\kern-1mm\char'022}} 
\def\iimage{\hbox{\font\=cmr7\\kern.3mm\char'022\kern-.7mm\char'022\kern-.3mm}} 
\def\images{\hbox{\font\=cmr10\\char'022\kern-1mm\char'022\kern-1mm\char'022}} 
\def\timesn{\kern-.2mm\times\kern-.2mm} 
\def\ttimes{\hbox{\kern-.2mm${}\times\kern-2.5mm\lower.8mm\hbox{\font\=cmr5\t}\kern1.8mm$}} 
\def\ttimesn{\hbox{\kern-.2mm${}\times\kern-2.5mm\lower.8mm\hbox{\font\=cmr5\t}\kern1.4mm$}} 
\def\ktimes{\hbox{\kern-.2mm${}\times\kern-2.5mm\lower1mm\hbox{\font\=cmr5\k}\kern1.5mm$}} 
\def\vstimes{\kern.95mm\raise.45mm\hbox{\font\=cmbsy6\\char'002}\kern-2.3mm\lower.9mm\hbox{\font\=cmr5\vs}\kern1.05mm} 

\def\Examplee{{\font\=cmssi10\E\kern.15mmx\kern.15mma\kern.15mmm\kern.14mmp\kern.17mml\kern.15mme}\kern.3mm. }
\def\Examples{{\font\=cmssi10\E\kern.15mmx\kern.15mma\kern.15mmm\kern.14mmp\kern.17mml\kern.15mme\kern.15mms}\kern.3mm. }
\def\Remarkk{{\font\=cmssi10\R\kern.15mme\kern.15mmm\kern.15mma\kern.15mmr\kern.15mmk\kern.15mm. }}
\def\Remarkss{{\font\=cmssi10\R\kern.15mme\kern.15mmm\kern.15mma\kern.15mmr\kern.15mmk\kern.15mms}\kern.3mm. }
\def\No{{I\!\!N\kern-.54mm\lower.15mm\hbox{$_{\rm o}$}}} 
\def\Nopot#1{I\!\!N\kern-.54mm\lower.15mm\hbox{$_{\rm o}$}\kern-.7mm{}^{#1}} 
\def\potNo{^{\kern.37mm I\!\!{N_{}}_{\kern-.22mm{\rm o}}}} 
\def\minus{\kern.2mm\lower1.05mm\hbox{$^-$}}
\def\pplus{\raise.22mm\hbox{\font\=cmr5\\char'053}}
\def\mminus{\raise.18mm\hbox{\font\=cmsy5\\char'000}}
\def\plusinftyy{\raise.18mm\hbox{\font\=cmr5\\char'053}\infty}
\def\minusinftyy{\raise.18mm\hbox{\font\=cmsy5\\char'000}\infty}
\def\plusinfty{\lower1.05mm\hbox{$^+$}\infty}
\def\minusinfty{\lower1.05mm\hbox{$^-$}\infty}
\def\Qe{\hbox{$Q\kern-2.6mm\raise.2mm\hbox{\font\=cmssqi8\I}\kern1.7mm$}}
\def\Re{{I\!\!R}}
\def\Rep{{{I\!\!R^{\phantom{l}}}^{{}_{{}^{\!}\!+}}}}
\def\Repp{{{{{{{I\!\!R_{}}_{}}_{}}_{}}}_{{}^{\!+}}}}
\def\Ce{{\hbox{$C\kern-2.5mm\raise.2mm\hbox{\font\=cmssqi8\I}\kern1.48mm$}}}
\def\imag{\kern.15mm\lower.6mm\hbox{$^{^*}$}\kern-1.8mm\imath\kern.1mm} 

\def\ebiF{\kern.1mm\hbox{\font\=cmmib8\F}\kern.5mm} 
\def\ebiT{\kern.1mm\hbox{\font\=cmmib8\T}\kern.6mm} 
\def\ebiU{\kern.1mm\hbox{\font\=cmmib8\U}\kern.5mm} 
\def\biit#1{\hbox{\font\=cmmib10\#1}} 

\def\fssi#1{\hbox{\font\=cmssi10\#1}\kern0.15mm} 
\def\smb#1{\hbox{\font\†=cmmi8\†#1\kern.3mm}} 
\def\eCal#1{\kern.1mm\hbox{\font\†=cmbsy8\†#1\kern.4mm}} 
\def\ecal#1{\kern.1mm\hbox{\font\†=cmsy8\†#1\kern.3mm}} 
\def\ncal#1{\kern.1mm\hbox{\font\†=cmsy9\†#1\kern.3mm}} 
\def\vcal#1{\kern-.1mm\vec{\kern.2mm\hbox{\font\†=cmsy7\†#1}\kern.3mm}} 

\def\idv{\hbox{\font\=cmr10\id}\kern.25mm\lower.8mm\hbox{\font\=cmr7\v}\kern.3mm} 
\def\seq#1{\langle#1\rangle}

\def\ymp{{}^{}{\Cal N}_o\,}
\def\SemiNor{\Cal S_{_N}\kern0.15mm}
\def\vecs{\upsilon\kern-0.3mm\lower.15mm\hbox{$_s$}\kern0.2mm} 
\def\vecss{\hbox{\font\=cmitt10\v}\kern-0.1mm\lower.15mm\hbox{$_s$}\kern0.2mm} 
\def\bnull#1{\hbox{\font\=cmssbx10\0}{}_{\font\=cmmi6\lower.15mm\hbox{\kern-.1mm\#1\kern.15mm}}} 
\def\bzero#1{\hbox{\font\=cmbx10\0}{}_{\font\=cmmi6\lower.15mm\hbox{\kern-.1mm\#1\kern.15mm}}} 
\def\dom{{{}^{}{\rm dom}\,{}_{{}^{}}}}
\def\domm{{}^{}{\rm dom}^{\kern.3mm\hbox{\font\=cmr6\2}}\,}
\def\domr#1{\roman{dom}^{\font\=cmr6\raise.0mm\hbox{\kern.3mm\#1}}}

\def\rng{{}^{}{\rm rng}\,{}_{{}^{}}}
\def\CPi#1{C\kern-.2mm\lower.05mm\hbox{$_{_\Pi}$}\kern-1.52mm{}^{#1}}
\def\CinftyPi{C\kern.4mm\raise.3mm\hbox{$^\infty$}\kern-3.35mm_{_\Pi}\kern1.45mm}
\def\CinftyS{\Cinfty\kern-3.9mm_{_{\Cal S}}\kern1.45mm}
\def\Cinfty{C\kern.4mm\raise.3mm\hbox{$^\infty$}\kern.15mm}
\def\Cinftyzero{\hbox{$C\kern.4mm\raise.3mm\hbox{$^\infty$}\kern.15mm\kern-3.5mm_{\font\=cmr6\lower.15mm\hbox{\kern.1mm\0}}\kern1.9mm$}}
\def\Lip#1{{}^{}{{{{{{\Cal L}{}^{}ip}^{}}^{}}^{}}^{}}^{#1{}^{}}}
\def\RHB#1#2{\raise#1mm\hbox{$#2$}} 
\def\LHB#1#2{\lower#1mm\hbox{$#2$}} 

\def\fiveroman#1{\hbox{\font\=cmr5\#1\kern.1mm}}

\def\sixroman#1{\hbox{\font\=cmr6\#1\kern.1mm}}
\def\eightmath#1{\hbox{\font\=cmmi8\{#1}\kern.1mm}}
\def\eightroman#1{\hbox{\font\=cmr8\{#1}\kern.1mm}}
\def\erm#1{{\font\=cmr8\#1}}
\def\subtext#1{\raise.2mm\hbox{$_{_{\kern0.15mm\roman{#1}}}$}}
\def\subtexT#1{\raise.2mm\hbox{$_{_{\kern0.15mm\hbox{\font\=cmr5\#1}}}$}}
\def\sNor#1{\kern.25mm\lower.38mm\hbox{$_{#1}$}}
\def\sNorr#1{\kern-.2mm\lower.38mm\hbox{$_{#1}$}}
\def\sNoreset_#1{\kern.13mm\lower.83mm\hbox{\font\=cmmi6\C}\kern.32mm\lower.1mm\hbox{$_{^{\emptyset,#1}}$}}
\def\sbi#1{{_{\kern-0.1mm}}_{#1}} 
\def\ais#1_#2{{}_{\font\=cmmi6\lower.15mm\hbox{\kern-.1mm\#1\kern.15mm}}\lower.3mm\hbox{${_{\kern-0.3mm_{#2}}}$}} %
\def\aais#1_#2{\kern.1mm{}_{\font\=cmmi6\lower.25mm\hbox{\kern-.1mm\#1\kern.15mm}}\lower.4mm\hbox{${_{\kern-0.3mm_{#2}}}$}} %
\def\ai#1{{}_{\font\=cmmi6\lower.15mm\hbox{\kern-.1mm\#1\kern.15mm}}} 
\def\yi#1{^{\font\=cmmi6\raise.0mm\hbox{\kern-.1mm\#1\kern.15mm}}} 
\def\ear#1{{}_{\font\=cmr5\lower.15mm\hbox{\kern.1mm\#1}}} 
\def\ar#1{{}_{\font\=cmr6\lower.15mm\hbox{\kern.1mm\#1}}} 
\def\aar#1{_{\font\=cmr6\lower.15mm\hbox{\kern.1mm\#1}}} 
\def\yr#1{^{\font\=cmr6\raise.0mm\hbox{\kern.3mm\#1}}} 
\def\yrai^#1_#2{^{\kern.4mm\hbox{\font\=cmr6\{#1}}}_{\kern.2mm{#2}}}
\def\upparentes#1{^{\kern.2mm\raise.2mm\hbox{\font\=cmr6\\char'050}\kern.1mm{#1}\kern.1mm\raise.2mm\hbox{\font\=cmr6\\char'051}}} 
\def\lupar{\kern.2mm\lower1mm\hbox{$^{^(}$}} 
\def\rupar{\lower1mm\hbox{$^{^)}$}\kern-.15mm} 
\def\yyi#1{^{\font\=cmmi6\lower.6mm\hbox{\kern-.25mm\#1\kern-.05mm}}} 
\def\yyr#1{^{\font\=cmr6\lower.45mm\hbox{\kern-.25mm\#1\kern-.15mm}}} 
\def\yplus{\lower1mm\hbox{$^{^+}$}} 
\def\yminus{\lower1mm\hbox{$^{^-}$}} 
\def\aminus{{\kern.15mm\raise.3mm\hbox{$_{_-}$}\kern-.1mm}}%
\def\adot{\kern.2mm\hbox{\font\=cmb10\\char'056}}%
\def\ydot{\kern.2mm\raise1.9mm\hbox{\font\=cmb10\\char'056}}
\def\yydot{\kern.2mm\raise1.35mm\hbox{\font\=cmb7\\char'056}\kern.2mm}
\def\yydott{\kern.2mm\raise1.35mm\hbox{\font\=cmb6\\char'056}\kern.2mm}
\def\ClT{{\rm Cl}\kern.25mm\lower.4mm\hbox{$_{\Cal T}$}\kern0.2mm} 
\def\IntT{\sp{\rm Int}\kern.2mm\lower.4mm\hbox{$_{\Cal T}$}\kern0.2mm} 
\def\Cl_taurd#1{\roman{Cl_{}}_{\kern0.37mm\hbox{\font\=cmmi8\\char'034}\kern-0.15mm{_{}}_{rd}\kern0.2mm#1\,}}
\def\inc{\subseteq}
\def\iinc{\supseteq}
\def\exi#1{\exists\,#1\kern.2mm\,;}
\def\all#1{\forall\,#1\kern.2mm\,;}

\def\equivv{\Leftrightarrow}

\def\spp{\kern0.07mm} 
\def\sp{\kern0.15mm} 
\def\ssp{\kern0.37mm} 
\def\snn{\kern-0.2mm} 
\def\sn{\kern-0.3mm} 
\def\ssn{\kern-0.63mm} 
\def\biggerlineskip#1 {\linebreak\nopagebreak\vskip-4.2mm\vskip.#1mm\nopagebreak\noindent}%
\def\Biggerlineskip#1 {\linebreak\nopagebreak\vskip-4.2mm\vskip#1pt\nopagebreak\noindent}%
\def\KP#1{\kern#1mm} 
\def\KN#1{\kern-#1mm} 
\def\nhskip#1mm{$\null$\kern#1mm}
\def\hyppy#1{$\phantom{}$\hskip#1}
\def\mhyppy#1{\null\kern#1mm}

\def\text#1{\hbox{\rm#1}}

\def\VBOX/#1/#2/HEREend{\vbox{#2\vskip-#1mm}\vfill\null\eject}
\def\œ$#1${\hbox{$#1$}} 
\def\"{\"a} \def\"{\"o}
\def\q#1{``\kern0.37mm#1\kern0.37mm"}
\def\newProCla#1\par#2\par{\vskip1.7mm\noindent\bf#1\it#2\vskip1.7mm}
\def\Prooff{{\font\=cmssi10\P\kern.37mmr\kern.37mmo\kern.37mmo\kern.37mmf\kern.37mm. }\rm}

\def\QED{\hfill\hbox{$\ \sqcap\kern-2.45mm\sqcup$}}

\def\noin{\noindent}
\def\Newline{\kern-10mm\newline}
\font\rp=cmr8
\def\eps{\varepsilon}
\def\leu{\raise1.5mm\hbox{\font\=cmmi5\\char'074}\kern.2mm}%
\def\riu{\kern.2mm\raise1.5mm\hbox{\font\=cmmi5\\char'076}}%
\def\Symbol#1Ï{\kern.35mm\hbox{\font\=cmr10\\char'047}\kern.2mm#1\kern.35mm\hbox{\font\=cmr10\\char'047}}%
\def\Symboo#1Ï{\kern.35mm\text{`}\kern.2mm#1\kern.35mm\hbox{\font\=cmr10\\char'047}}
\def\RunMyHead#1#2#3#4{%
 \headline{\ifnum\pageno=\firstpage\hfil%
           \else{\ifodd\pageno{\rp#3\phantom\folio\hfil#4\hfil\phantom{#3}\folio}%
                 \else{\rp\folio\phantom{#2}\hfil#1\hfil\phantom\folio#2}%
                 \fi}%
           \fi}%
 \footline{\ifnum\pageno=\firstpage\hfil{\rp[\,\folio\,]}\hfil%
           \else\hfil%
           \fi}%
}%
\def\bulgin{\noindent$\bullet$ \ \kern.1mm} 
\def\bulgen{\noindent\kern-1.5mm$\bullet$\kern3.95mm} 
\def\subhead#1\par#2\par{\vskip4mm\smallbreak\null\smallskip\vbox{\noindent\bbf#1\hfill\kern1.5mm#2\hfill\phantom{#1}\vskip2.5mm\nopagebreak}\nopagebreak\noindent}
\def\subheadd#1\par#2\par#3\par{\vskip4mm\smallbreak\null\smallskip\vbox{\noindent\bbf#1\hfill#2\hfill\phantom{#1}\vskip1.5mm\centerline{#3}\vskip2.5mm\nopagebreak}\nopagebreak\noindent}
\def\insubsubhead#1\par{\vskip4mm$\null$\hskip2mm{\font\=cmss10\#1}\vskip2mm\noindent}%
\def\binsubsubhead#1#2\par{\vskip4mm{\bf#1.}\hskip5mm{\font\=cmss10\#2}\vskip2mm\noindent}%
\def\wave{\hbox{\font\†=cmsy10\†\hbox{\char'164}\kern-2.35mm\hbox{\char'165}\kern.55mm}}
\def\wavee{\hbox{\font\†=cmsy8\†\hbox{\char'164}\kern-2.0mm\hbox{\char'165}\kern.4mm}} 
\def\barmj{\kern.25mm\bar{\hbox{\font\=cmr10\\char'021}}\kern.4mm}

\def\sigrd{\sigma\kern-.3mm_{_{rd}}\kern.15mm} 
\def\ssigrd{\sigma\kern0.45mm^{\font\=cmr6\hbox{\2}}\kern-2.2mm_{_{rd}}\kern.15mm}
\def\sssigrd{\sigma\kern0.45mm^{\font\=cmr6\hbox{\3}}\kern-2.2mm_{_{rd}}\kern.15mm}
\def\taurd{\tau_{_{rd}}\kern.15mm}
\def\tsigrd{\tau\sigma\kern-.3mm_{_{rd}}\kern.15mm} 
\def\tauR#1{\tau_{_{I\!\!R}}\kern-1.5mm^{#1}}
\def\RN{I\!\!R\kern.3mm^{\hbox{\font\=cmmi6\N}}} 
\def\QTN{Q\kern.1mm_{\lower.2mm\hbox{\font\=cmmi6\T}}^{\kern.2mm\hbox{\font\=cmmi6\N}}} 

\def\leLCS-{{\le}{}_{_{{\rm LCS}}}\text{\sp-\sp}}
\def\Centerline#1\par#2\par#3{\noindent#1\phantom{#3}\hfill#2\hfill\phantom{#1}#3}

\def\C_c^#1(#2){C_{\roman c}^{\kern.6mm#1}\kern.2mm(\boldsymbol#2\kern.4mm)}

\def\Lip_#1t^#2{\Cal Lip_{\,\fiveroman{#1\ssp t}}^{\kern.9mm#2}}
\def\variat_#1^#2 {\hbox{\font\=cmtex11\\char'012}\sbi{#1^{\vphantom l}\,}^{\kern1mm#2}}
\def\variamap^#1{\hbox{\font\=cmtex11\\char'012}^{\kern.7mm#1}\kern-.2mm}
\def\tauMac{\tau\kern.15mm\lower.7mm\hbox{\font\=cmr5\Mac}\kern.6mm}
\def\cdotn{\snn\cdot\snn}
\def\uprime{\kern.4mm\lower.7mm\hbox{$^{'}$}}

\newcommand\bosy{\boldsymbol}

\begin{document}

\title[$         \text{\sc Inverse function theorem}$]%
                    {An inverse function theorem for\vskip1mm
             Colombeau tame Fr\"licher\,--\,Kriegl maps}

\author[S. Hiltunen]{Seppo\ I\. Hiltunen}
\address{Helsinki University of Technology                             \vskip0mm$\hspace{2mm}$
           Institute of Mathematics, U311                              \vskip0mm$\hspace{2mm}$
           P.O.\ Box 1100                                              \vskip0mm$\hspace{2mm}$
           FIN-02015 HUT\vskip0mm
         FINLAND}
\email{shiltune\,@\,cc.hut.fi}



\begin{abstract}

For \œ$k=1\sp,2\ssp,\ldots\,\infty$ and a Fr\"licher\,--\,Kriegl order $k$
Lipschitz differentiable map \œ$f:E\iinc U\to E$ having derivative at \œ$
x\ear 0\in U$ a linear homeomorphism \œ$E\to E$ and satisfying a Colombeau
type tameness condition, we prove that $x\ear 0$ has a neighbourhood \œ$V\inc
U$ with $f\,|\,V$ a local order $k$ Lipschitz diffeomorphism. As a corollary
we obtain a similar result for Keller $C_c^{\,\infty}$ maps with $E$ in a
class including Fr\'echet and Silva spaces. We also indicate a procedure for
verifying the tameness condition for maps of the type
\œ$x\mapsto\varphi\circ[\,\roman{id\,},x\,]$
and spaces \œ$E=C^{\,\infty}(Q)$ when $Q$ is compact by considering the case \œ$
Q=[\,0\,,1\,]\,$. Our considerations are motivated by the wish to try to
retain something valuable in an interesting but defective treatment of
integrability of Lie algebras by J. Leslie.                   \end{abstract}

\maketitle


\noin In \cite[Theorem 4.1, p.\ 439]{Les} an interesting assertion is given
which via \cite[Lemma 4.2, p.\ 441]{Les} is based on
\cite[Theorem 2.2, p.\ 430]{Les} whose proof in turn is omitted for the most
part with the exception of just few hints. It even seems that it cannot be
proved unless one there requires the space $E$ to be suitably restricted. To
fill at least this gap in \cite{Les}\ssp, we prove Theorem \ref{inve func}
below which then gives Corollary \ref{Frec inve fu} as a replacement for
\cite[Theorem 2.2]{Les} when $E$ is cm-convenient.

Proposition \ref{appl exa} below should serve as a prototype for assertions
guaranteeing that Corollary \ref{Frec inve fu} can be applied to maps of the
type \œ$\,x\mapsto\varphi\circ[\,\sp\roman{id\,},x\,]\ssp$ of spaces $\ssp
\Cinfty(Q)$ when $Q$ is compact. Below, we shall use the notational
conventions of \cite{Hic} from which we in particular recall the following

\begin{conventions}

Letting $\biit R$ be the topological field of real numbers, the class of all\linebreak
real locally convex Hausdorff topological vector spaces is $\LCS(\biit R\sp)\,
$. For \œ$E\in\LCS(\biit R\sp)$ we have \œ$E=(X\sp,\Cal T\ssp)$ where \œ$X =
\sigrd E$ is the underlying \q{abstract} real vector space, and \œ$\sp\Cal T =
\taurd E$ is the (locally convex) topology for the underlying set $\vecs E\,$.
The filter of $\,\Cal T\,$--\,neighborhoods of the zero vector $\bnull E$ is $
\ymp E\,$, and the von Neu- mann bornology, the set of all bounded sets is $
\rajou E\,$.

A vector {\it map\ssp} of spaces in $\LCS(\biit R\sp)$ is any triplet \œ$
\tilde f = (E\ssp,F\sp,f\sp)$ such that \œ$E\ssp,F\in$ $\LCS(\biit R\sp)$ and
$f$ is a function with \œ$f\inc(\sp\vecs E\sp)\timesn(\sp\vecs F\sp)\,$. The
function value of $f$ at $x$ is $f\sp\fvalue x\ssp$, instead of the
conventional notation \q{\snn$f\sp(x)$}. The derivative at $x$ of a Gateaux
differentiable map $\tilde f=(E\ssp,F\sp,f\sp)$ is $\tilde f\ssp'(x)\,$.

By definition, we have $\domm\Gamma=\dom\snn(\dom\Gamma\sp)$ for any class $\ssp
\Gamma\ssp$.                                               \end{conventions}

For \œ$k\in\infty\ssp\yplus = \No\cup\{\infty\}\,$, we below consider the
differentiability classes $\Lip_FKt^k$ and $\C_c^k(R)\,$. The latter one has
as its members exactly the maps $(E\ssp,F\sp,f\sp)$ with \œ$E\ssp,F \in
\LCS(\biit R\sp)$ and \œ$\dom\sn f\in\taurd E$ such that for \œ$l\in k+1\adot$
the order $l$ variation $\variat_E\sp F^l f$ has \œ$
\dom\snn(\sp\variat_E\sp F^l f\sp) = (\sp\vecs E\ssp)^{\,\sp l\ssp+\sp 1.}
\cap \{\,\biit x:\biit x\fvalue\emptyset\in\dom\sn f\sp\,\}\,$, and $
\variat_E\sp F^l f$ is contin- uous $\taurd(E\expnota^l\ssp+\sp 1.]_{tvs})\to
\taurd F\sp$. For more information, see \cite[Section 3]{Hic} and \cite{Ke}\ssp.

The class $\Lip_FKt^k$ has as its elements exactly the maps $(E\ssp,F\sp,f\sp)$
such that the spaces $E\ssp,F\in\LCS(\biit R\sp)$ are bornological and
locally (i.e.\ Mackey) complete, see \cite[p.\ 196\,]{Jr} or
\cite[Lemma 2.2, p.\ 15\,]{KM}\ssp, and \œ$f:E\iinc\dom\sn f\to F$ is $
\Cal L\spp ip\,^k$ in the sense of \cite[pp.\ 83, 99]{FK} or
\cite[Definition 12.1, p.\ 118]{KM}\ssp. Our choice corresponds to the one
made in \cite{FK} where the spaces are bornological (locally convex) which is
not required in \cite{KM}\ssp. This has the consequence that for any fixed \œ$
k\in\infty\ssp\yplus$ the class $
\Lip_XZt^k\KN2\image\{\ssp(\biit R\,,E\ssp)\ssp\}$ of local $
k^{\,\fiveroman{th}}$ order Lipschitz differentiable curves in $E$ uniquely
determines $E$ in the class \œ$\LCS(\biit R\sp)\cap\{\,F:\sigrd E=\sigrd F\text{
and }\rajou E=\rajou F\sp\,\}$ when $\sbi{\,\fiveroman{XZ}\ssp =
\ssp\fiveroman{FK}}\,$, but not when $\sbi{\,\fiveroman{XZ}\ssp=
\ssp\fiveroman{KM}}\,$.

For \œ$E\in\LCS(\biit R\sp)\,$, a set \œ$U\inc\vecs E$ we call {\it mopen\ssp}
in $E$ if and only if for all \œ$x\in U$ and \œ$B\in\rajou E$ there is \œ$
\delta\in\Rep$ with \œ$t\,B\inc U-x$ for \œ$0\le t\le\delta\ssp$. The set \œ$
\tauMac E=$ \œ$\{\,U:U\text{ mopen in }E\sp\,\}$ then is a topology for $
\vecs E$ with \œ$\taurd E\inc\tauMac E\,$, equality here holding if $E$ is
metrizable, cf.\ \cite[Remark 2.4.5, p.\ 52]{FK} where $\tauMac E$ is called
the {\it Mackey closure topology\ssp}. For \œ$(E\ssp,F\sp,f\sp)\in\Lip_FKt^k$
it holds that \œ$\dom\sn f\in\tauMac E$ and that $f$ is continuous \œ$
\tauMac E\to\tauMac F\sp$, and further \œ$(E\ssp,F\sp,f\,|\,U\ssp) \in
\Lip_FKt^k$ for any \œ$U\in\tauMac E\,$. For these facts, we refer the reader
to see \cite[Proposition 2.3.7, p.\ 44, Corollary 4.1.7, p.\ 85,
Proposition 4.3.2, p.\ 99]{FK}\ssp.

For short, a space \œ$E\in\LCS(\biit R\ssp)$ we call {\it cm-convenient\ssp}
if{}f also $E$ is bornological and locally complete, and such that \œ$
\tauMac(\spp E\expnota^k]_{tvs}) = \taurd\spp(\spp E\expnota^k]_{tvs})$ holds
for \œ$k\in\No\ssp$. By the above, all Fr\'echet spaces are cm-convenient, and
by \cite[Theorem 7.3.2\sp(1)\ssp, p.\ 99]{HN2} also Silva spaces are
cm-convenient, cf.\ \cite[Theorem 6.1.4, p.\ 190]{FK}\ssp.

\begin{proposition}\label{rela wtw Lip^infty and C_c^infty}

Let \œ$\sp\tilde f=(E\ssp,F\sp,f\sp)$ where \œ$\sp F\in\domm\Lip_FKt^{0.}$ and
$\,E$ is cm-conven- ient. The equivalence $\,\tilde f\in\Lip_FKt^\infty\sp
\equivv\sp\tilde f\in\C_c^\infty(R)\,$ then holds.         \end{proposition}

\begin{proof} First letting \œ$\sp\tilde f\in\Lip_FKt^\infty\,$, for \œ$k \in
\No$ by \cite[Theorem 4.3.27, p.\ 112]{FK}\ssp, for the variation map \œ$
\variamap^k\tilde f = (E\expnota^k\ssp+\sp 1.]_{tvs}\sp,F\sp,
\variat_E\sp F^k f\ssp)$ we have \œ$\variamap^k\tilde f\in\Lip_FKt^{0.}$ with
\œ$\dom\snn(\sp\variat_E\sp F^k f\sp) = (\sp\vecs E\sp)^{\,k\ssp+\sp 1.} \cap
\{\,\biit x:\biit x\fvalue\emptyset\in\dom\sn f\sp\,\}\,$. By
\cite[Corollary 4.1.7, p.\ 85]{FK}\ssp, hence $\variat_E\sp F^k f$ is
continuous \œ$\tauMac(E\expnota^k\ssp+\sp 1.]_{tvs})\to\tauMac F\sp$, whence
by \œ$\taurd F\inc\tauMac F\sp$, also continuous \œ$
\taurd\spp(E\expnota^k\ssp+\sp 1.]_{tvs})\to\taurd F\sp$. So we get \œ$
\variamap^k\tilde f\in\C_c^{0.}(R)\,$. Here \œ$k\in\No$ being arbitrary, we
get \œ$\tilde f\in\C_c^\infty(R)\,$. Conversely, then letting \œ$\tilde f \in
\C_c^\infty(R)\,$, for an arbitrary \œ$
c\in\Lip_FKt^\infty\KN1\image\{\ssp(\biit R\,,E\ssp)\ssp\} =
         \C_c^\infty(R)\image\{\ssp(\biit R\,,E\ssp)\ssp\}$ the chain rule
gives the result that \œ$f\circ\sp c \in
\C_c^\infty(R)\image\{\ssp(\biit R\,,F\ssp)\ssp\}\,$. Directly by definition,
this further gives $\sp\tilde f\in\Lip_FKt^\infty\,$, in view of
\cite[Lemma 4.3.1, p.\ 99]{FK}\ssp.                              \end{proof}

\begin{definitions}

For \œ$E\in\LCS(\biit R\sp)\,$, a set \œ$\Cal B \inc
\rajou E\setminus\sn\{\emptyset\}$ of absolutely convex closed sets we call a
{\it bornological generator\ssp} for $E$ if{}f for all \œ$A\,,B\in\Cal B$
there is \œ$C\in\Cal B$ with $A\cup B\inc C\ssp$, and every $A\in\rajou E$ has
some $t\in\Rep$ and $B\in\Cal B$ with $A\inc t\,B\ssp.$

For any real vector space $X$ and any nonempty absolutely convex $B$ in $X$
and any \œ$x\in\vecss X$ with \œ$\,\roman S\sbi{X\sp B}\,x = \Rep\cap\{\,t :
t^{\sp\mminus 1}\sp x\in B\sp\,\}$ we let \œ$\|\ssp x\ssp\|\sNor{X\sp B} =
\inf\ssp(\sp\roman S\sbi{X\sp B}\,x\sp)\,$, hence having \œ$
\|\ssp x\ssp\|\sNor{X\sp B} = \plusinfty$ in case \œ$x\not\in S$ for the
linear span \œ$S = \taurd\sp X\sp\,[\,\Rep\sn\timesn B\sp\,]$ \œ$ = \{\,t\,v:
t\in\Rep\text{ and }v\in B\sp\,\}\,$. In particular, if we have \œ$X=\sigrd E$
with $E$ locally complete and also $B$ is $\taurd E\,$--\,closed with \œ$B \in
\rajou E\,$, by \cite[Proposition 10.2.1, p.\ 197]{Jr} then \œ$
(X_{\sp|\ssp S}\ssp,\seq{\,\|\ssp x\ssp\|\sNor{X\sp B}\sn:x\in S\sp\,}\spp)$
is a normed Banach space with $B$ its closed unit ball. The corresponding {\it
Banachable\ssp} locally convex topological vector space is \œ$X\sbi B =
(X_{\sp|\ssp S}\ssp,\Cal T\ssp)$ for \œ$\sp\Cal T=\Pows S\cap\{\,U :
\all{x\in U}\,\exi{\eps\in\Rep}\,\eps\,B\inc U-x\,\}\,$. Note that we have $\,
\tauMac E\sp\lei S\inc\Cal T\spp$.                         \end{definitions}

The standard argument in \cite[Theorem 10.7, pp.\ 231\,--\,232]{Rud:FA} gives
the following

\begin{lemma}\label{Neuma ser}

Let \œ$\ssp 0\le\eps<1\ssp,$ and let $\sp(X\sp,\Nu\sp)$ be a normed Banach
space. If also $\,\ell\sp$ is linear \œ$X\to X$ with \œ$\sp
\Nu\sp\fvalue(\ssp\ell\ssp\fvalue x-x\sp)\le\eps\,(\sp\Nu\sp\fvalue x\sp)$ for
all \œ$\ssp x\in\vecss X\sp,$ then $\,\ell\ssp\inve$ is linear $\sp X\to X\sp$
with $\mhyppy9
\Nu\circ\ell\ssp\inve\fvalue x \le
 (\sp 1-\eps\sp)^{\sp\mminus 1}\sp(\sp\Nu\sp\fvalue x\sp)\mhyppy6$ \ and \newline $\mhyppy{21.8}
\Nu\sp\fvalue(\ssp\ell\ssp\inve\fvalue x-x\sp) \le
 (\sp 1-\eps\sp)^{\sp\mminus 1}\ssp\eps\,(\sp\Nu\sp\fvalue x\sp)\mhyppy5$
for all $\,x\in\vecss X\sp$.                                     \end{lemma}

\begin{definitions}\label{Colo 0}

We let Colo$\ar A\sp(\sp\tilde f\sp,y\ar 0\ssp,\eps\ssp,\Cal B\sp)$ mean that
there are $E\ssp,f$ such that \œ$\tilde f=$ $(E\ssp,E\ssp,f\sp)$ with $
\tilde f$ a Gateaux differentiable map and \œ$E\in\LCS(\biit R\sp)$
bornological and locally complete and \œ$\eps\in\Rep$ and $\Cal B$ a
bornological generator for $E$ such that\linebreak for \œ$v\in B\in\Cal B$ and
for \œ$y\ar 1 = y\ar 0+2\,v$ we have \œ$y\ar 1 \in (\dom\sn f\sp)\cap(\sp
\vecs\spp(\sp\sigrd E\sbi B\spp))\,$, and also the inclusion $\,
    (\sp\tilde f\ssp'(\sp y\ar 1)-\tilde f\ssp'(\sp y\ar 0))\image B \inc
               \eps\,\tilde f\ssp'(\sp y\ar 0)\image B\,$ holds.

Let Colo$\ar{\,0}\ssp(\sp\tilde f\sp,y\ar 0\ssp,\eps\ssp,\Cal B\sp)$ mean
that Colo$\ar A\sp(\sp\tilde f\sp,y\ar 0\ssp,\eps\ssp,\Cal B\sp)$ holds and
$\sp\tilde f\ssp'(\sp y\ar 0)$ is a linear homeomorphism \œ$E\to E\,$. For \œ$\sp
\tilde f=(E\ssp,E\ssp,f\sp)\,$, letting \œ$\sp\tilde f\,|\subtext{map\,}V =
(E\ssp,E\ssp,f\,|\,V\ssp)\,$, a vector map $\tilde f$ we say to be {\it
Colombeau$\ar{\,0\,}$--\,tame\ssp} at $y\ar 0$ if{}f for all $\,V$ with \œ$
y\ar 0\in V\sn\in\taurd E\ssp$ there is some $\Cal B$ such that \ssp
Colo$\ar{\,0}\ssp(\ssp\tilde f\,|\subtext{map\,}V\spp,\sp y\ar 0\,,\sp
\frac12\,,\Cal B\ssp)\ssp$ holds.                          \end{definitions}

\begin{remark}\label{notes abou Colo 0}

We note some facts to be used below.

(a) \ Assuming \ssp Colo$\ar{\,A}\ssp(\sp\tilde f\sp,y\ar 0\ssp,\eps\ssp,
    \Cal B\sp)$ to hold with \œ$\tilde f=(E\ssp,E\ssp,f\sp)$ and \œ$X =
\sigrd E$\linebreak and \œ$B\in\Cal B\ssp$, since \œ$\bnull E\in B$ by
circledness of $B$ following from its absolute convexity, we have \œ$\sp
y\ar 0 = y\ar 0+2\,\bnull E\in\dom\sn f\spp\cap\sp\vecs(X\sbi B\ssp)\,$, and
consequently \œ$Q\inc\vecs(X\sbi B\ssp)$ holds for the closed convex set \œ$
Q = y\ar 0+2\,B\ssp$. Further \œ$(\sp\taurd\spp(X\sbi B\spp)\ssp,\taurd E\ssp,
f\,|\,Q\sp)$ is a {\it topological map\,}, i.e.\ we have $\sp f\,|\,Q$
continuous $\taurd\spp(X\sbi B\spp)\to\taurd E\,$.

To see this, arbitrarily fixing \œ$y\ar 1\in Q$ and a closed convex \œ$V \in
\ymp E\,$, there should be some \œ$N\in\Nbh(\sp y\ar 1\sp,
\taurd\spp(X\sbi B\ssp)\lei Q\sp)$ with \œ$\sp f\sp\image\snn N \inc
f\sp\fvalue y\ar 1\snn+V\spp$. For this, we first note that $
\tilde f\ssp'(\sp y\ar 0)$ being a continous linear map \œ$E\to E\,$, it is
bornological \œ$\rajou E\to\rajou E\,$, and hence we have \œ$
\tilde f\ssp'(\sp y\ar 0)\image\snn B\in\rajou E\,$, whence further there is
some \œ$\eps\ar 1\in\Rep$ with \œ$
\eps\ar 1\ssp\tilde f\ssp'(\sp y\ar 0)\image\snn B\inc V\spp$. With \œ$
\lambda = (\sp 1+\eps\sp)^{\sp\mminus 1}\ssp\eps\ar 1\ssp$, now taking \œ$ N =
(\sp y\ar 1\snn+\lambda\,B\sp)\cap Q\ssp$, for\linebreak $y\ar 2 =
y\ar 1\snn+v\in N$ we see $\sp f\sp\fvalue y\ar 2 \in
               f\sp\fvalue y\ar 1\snn+V$ to hold as follows.

For \œ$I=[\,0\,,1\,]$ and \œ$\,c\sp =
\seq{\,f\sp\fvalue(\sp y\ar 1\snn+t\,v\sp) -
       f\sp\fvalue y\ar 1\sn:t\in I\sp\,}\,$, we have $c$ a differentiable
curve in $E$ with \œ$\sp c\fvalue 0=\bnull E$ and \œ$\,c\fvalue 1 =
f\sp\fvalue y\ar 2 - f\sp\fvalue y\ar 1\ssp$. By the mean value theorem, it
hence suffices that $\rng\snn(\sp\roman D\sbi E\,c\ssp)\inc V\spp$. This is
the case since for $t\in I$ we have\vskip.5mm$\mhyppy{-2}
\roman D\sbi E\,c\fvalue t = \tilde f\ssp'(\sp y\ar 1\snn+t\,v\sp)\fvalue v
 \in (\sp\tilde f\ssp'(\sp y\ar 1\snn+t\,v\sp) -
\tilde f\ssp'(\sp y\ar 0))\image(\sp\lambda\,B\sp) +
\tilde f\ssp'(\sp y\ar 0)\image(\sp\lambda\,B\sp)  $\vskip.3mm\centerline{$\phantom{\roman D\sbi E\,c\fvalue t}
 \inc \lambda\,(\sp\eps\,\tilde f\ssp'(\sp y\ar 0)\image\snn B +
\tilde f\ssp'(\sp y\ar 0)\image\snn B\sp) \inc
\lambda\,(\sp 1+\eps\sp)\,\tilde f\ssp'(\sp y\ar 0)\image\snn B
\inc \eps\ar 1\ssp\tilde f\ssp'(\sp y\ar 0)\image\snn B \inc V\spp$.}\vskip.5mm

(b) \ Assuming \ssp Colo$\ar{\,0\,}(\sp\tilde f\sp,y\ar 0\ssp,\eps\ssp,
    \Cal B\sp)$ to hold with \œ$\tilde f=(E\ssp,E\ssp,f\sp)$ and \œ$X =
\sigrd E$\linebreak and \œ$\sp\eps<1\ssp$, for \œ$\sp y\ar 1\in y\ar 0 +
2\,\bigcup\ssp\Cal B$ we have $\tilde f\ssp'(\sp y\ar 1)$ a {\it linear
homeomorphism\ssp} \œ$E\to E\,$. For this, since we assume $E$ to be
bornological, it suffices that $\tilde f\ssp'(\sp y\ar 1)$ is bijective $
\vecs E\to\vecs E\,$, and that $\tilde f\ssp'(\sp y\ar 1)$ and $
(\sp\tilde f\ssp'(\sp y\ar 1))\inve$ are bornological $\rajou E\to\rajou E\,$.

First, to show indirectly that $\tilde f\ssp'(\sp y\ar 1)$ is injective, if
not, there is \œ$v\in\vecs E\setminus\sn\{\ssp\bnull E\}$ with \œ$
\tilde f\ssp'(\sp y\ar 1)\fvalue v = \bnull E\,$. There further is \œ$B\ar 1
\in\Cal B$ with \œ$y\ar 1\in y\ar 0 + 2\,B\ar 1\ssp$, and there is \œ$B\ar 2
\in\Cal B$ with \œ$v\in\vecs(X\sbi{B_2})\,$. We then find \œ$B\ar 3\in\Cal B$
with \œ$B\ar 1\snn\cup B\ar 2\inc B\ar 3\,$, and we have \œ$\,
    (\sp\tilde f\ssp'(\sp y\ar 1)-\tilde f\ssp'(\sp y\ar 0))\image B\ar 3 \inc
               \eps\,\tilde f\ssp'(\sp y\ar 0)\image B\ar 3\,$.
For \œ$\ssp\ell\sp = (\sp\tilde f\ssp'(\sp y\ar 0))\inve\snn \circ
(\sp\tilde f\ssp'(\sp y\ar 1))\,|\,\vecs(X\sbi{B_3})$ it follows that \ $
0 < \|\ssp v\ssp\|\sNor{X\sp B_3} =
\|\,\sp\ell\ssp\fvalue v-v\,\|\sNor{X\sp B_3} \le
\eps\,\|\ssp v\ssp\|\sNor{X\sp B_3} < \|\ssp v\ssp\|\sNor{X\sp B_3}\,$.

To get \œ$\sp\vecs E\inc\rng\snn(\sp\tilde f\ssp'(\sp y\ar 1))\,$, given \œ$
v\in\vecs E\,$, deducing as above, we find $B\ar 3$ with \œ$\,
(\sp\tilde f\ssp'(\sp y\ar 1)-\tilde f\ssp'(\sp y\ar 0))\image B\ar 3 \inc
               \eps\,\tilde f\ssp'(\sp y\ar 0)\image B\ar 3\,$, and now \œ$
\tilde f\ssp'(\sp y\ar 0)\inve\fvalue v\in\vecs(X\sbi{B_2})\,$. For \œ$x \in
\vecs(X\sbi{B_3})$ then \œ$\|\,\sp\ell\ssp\fvalue x-x\,\|\sNor{X\sp B_3} \le
\eps\,\|\ssp x\ssp\|\sNor{X\sp B_3}$ whence Lemma \ref{Neuma ser} gives $
\tilde f\ssp'(\sp y\ar 0)\inve\fvalue v$ $\in \vecs(X\sbi{B_3}) \inc
\rng\sp\ell\sp\,$, \sp and so there is $u\in\vecs(X\sbi{B_3})\inc\vecs E$ with
$\tilde f\ssp'(\sp y\ar 1)\fvalue u=v\ssp$.

To prove that $\tilde f\ssp'(\sp y\ar 1)$ is bornological \œ$\rajou E\to\rajou
E\,$, given \œ$B\in\rajou E\,$, as above, we find $B\ar 3$ with now also
having \œ$B\inc\lambda\,B\ar 3$ for a suitable \œ$\lambda\in\Rep\sp$. Then we
obtain \œ$\tilde f\ssp'(\sp y\ar 1)\image\snn B \inc
\lambda\,\tilde f\ssp'(\sp y\ar 1)\image\snn B\ar 3 \inc
\lambda\,(\sp 1+\eps\sp)\,\tilde f\ssp'(\sp y\ar 0)\image B\ar 3\in\rajou E\,
$. To get the assertion for $(\sp\tilde f\ssp'(\sp y\ar 1))\inve\sp$, we
arrange \œ$B\inc\lambda\,\tilde f\ssp'(\sp y\ar 0)\image\snn B\ar 3\,$, and by
Lemma \ref{Neuma ser} we obtain \œ$\ell\ssp\inve\image\snn B\ar 3\inc$ $
(\sp 1-\eps\sp)^{\sp\mminus 1}\ssp B\ar 3\,$, whence $\ssp
\tilde f\ssp'(\sp y\ar 1)\inve\image\snn B \inc
\lambda\,(\sp 1-\eps\sp)^{\sp\mminus 1}\ssp B\ar 3\in\rajou E\,$ follows. \vskip.3mm

(c) \ Assuming \,Colo$\ar{\,0\,}(\sp\tilde f\sp,y\ar 0\ssp,\eps\ssp,
    \Cal B\sp)\,$ to hold with \œ$\,\tilde f=(E\ssp,E\ssp,f\sp)\,$ and \œ$\,
X = \sigrd E$\linebreak and \œ$\,\,\ell\sp =
(\sp\tilde f\ssp'(\sp y\ar 0))\inve\sp$, for \œ$\sp f\aar 1 =
\seq{\,y-\ell\sp\circ f\sp\fvalue y:y\in\dom\sn f\sp\,}\,$, and for \œ$
y\ar 8\ssp,y\ar 9\in y\ar 0+2\,B$ with \œ$B\in\Cal B\ssp$, we have \œ$
\|\,f\aar 1\KN1\fvalue y\ar 8 - f\aar 1\KN1\fvalue y\ar 9\ssp\|\sNor{X\sp B}
\le \eps\,\|\,y\ar 8 - y\ar 9\ssp\|\sNor{X\sp B}\,$. Indeed, putting\linebreak
\œ$v = y\ar 8 - y\ar 9$ and \œ$I=[\,0\,,1\,]\,$, and considering in the space
$E$ the differentiable\linebreak curve \œ$
c \sp = \seq{\,t\,v - \ell\ssp\fvalue(\sp f\sp\fvalue(\sp y\ar 9+t\,v\sp) -
f\sp\fvalue y\ar 9):t\in I\sp\,}\,$, for which we have \œ$c\fvalue 0=\bnull E$\linebreak
and \œ$c\fvalue 1 = f\aar 1\KN1\fvalue y\ar 8 - f\aar 1\KN1\fvalue y\ar 9$ and
\œ$\roman D\sbi E\,c = \seq{\,v-\ell\sp\circ(\sp\tilde f\ssp'(\sp y\ar 9 +
t\,v\sp))\fvalue v:t\in I\sp\,}\,$, to get\linebreak the assertion, by the
mean value theorem, for arbitrarily fixed \œ$t\in I\sp$, it suffices\linebreak
that \œ$\|\,\sp\roman D\sbi E\,c\fvalue t\,\|\sNor{X\sp B} \le
\eps\,\|\ssp v\ssp\|\sNor{X\sp B}\,$. This is the case, since for \œ$y\ar 7 =
y\ar 9+t\,v$ and for \œ$\,\ell\ar 1 =
\ell\sp\circ(\sp\tilde f\ssp'(\sp y\ar 7))\,$, we have the inclusion \œ$
(\sp\tilde f\ssp'(\sp y\ar 7)-\tilde f\ssp'(\sp y\ar 0))\image B \inc
           \eps\,\tilde f\ssp'(\sp y\ar 0)\image B\ssp$, which further gives
$\, \|\,\sp\roman D\sbi E\,c\fvalue t\,\|\sNor{X\sp B}
=   \|\,\sp\ell\ar 1\KN1\fvalue v-v\,\|\sNor{X\sp B}
\le \eps\,\|\ssp v\ssp\|\sNor{X\sp B}\,$.                       \end{remark}

\begin{lemma}\label{on impl fus}

Let {\,\rm Colo}$\ar{\,0}\ssp(\sp\tilde f\sp,y\ar 0\ssp,\eps\ssp,\Cal B\sp)$
hold with \œ$\ssp 0\le\eps\le\frac 12$ and \œ$\,\tilde f=(E\ssp,E\ssp,f\sp)$
and \œ$\,X=\sigrd E$ and $\,x\ar 0=f\sp\fvalue y\ar 0\,$.
Also let \œ$\,\ell\sp = (\sp\tilde f\ssp'(\sp y\ar 0))\inve$ and \œ$\,Q =
x\ar 0+\sp\ell\ssp\inve\image\bigcup\ssp\Cal B\,$. Then there is a function \œ$\,
g\inc f\inve$ with \œ$\sp x\ar 0\in\dom g\in\taurd E\sp\cap\Pows Q\ssp,$ and
in addition for every $\ssp B\in\Cal B$ and for all $\,x\ar 1\ssp,x\ar 2 \in
(\sp x\ar 0+\sp\ell\ssp\inve\image\snn B\ssp)\spp\cap\dom g\,$ it holds that\vskip.3mm\centerline{$
\|\,g\fvalue x\ar 1\snn-g\fvalue x\ar 2\ssp\|\sNor{X\sp B} \le
(\sp 1-\eps\sp)^{\sp\mminus 1}\sp\|\,\sp\ell\ssp
\fvalue(\sp x\ar 1\snn-x\ar 2)\,\|\sNor{X\sp B}\,$.}             \end{lemma}

\begin{proof} Let \ \ \ $
h=\{\ssp(\sp x\ssp,y\ssp,z\sp):\exi{x\ar 1}\,(\sp y\ssp,x\ar 1)\in f$ and $
\ell\ssp\fvalue x-\ell\ssp\fvalue x\ar 1\snn+y=z\,\}\,$.

\noin Then $h$ is a function, and we further put \œ$g\ar 1=\domm\Gamma\sp$,
where $\Gamma$ is the set of all $(\sp x\ssp,y\ssp,B\ssp,\biit y\sp)$ such
that \œ$B\in\Cal B$ and \œ$\biit y\in(\sp\vecs E\ssp)\potNo$ with \œ$x \in
x\ar 0+\tilde f\ssp'(\sp y\ar 0)\image B$ and\linebreak
$\biit y\sp\fvalue\emptyset=y\ar 0$ and $\biit y\to y$ in top $\taurd E$ and $
(\sp x\ssp,\biit y\sp\fvalue i\ssp,\biit y\sp\fvalue i\ssp\yplus)\in h$ for
all $i\in\No\ssp$.

We note that $g\ar 1$ is a function, since if \œ$(\sp x\ssp,y\sbi\iota) \in
g\ar 1$ for $\sbi{\iota\ssp=\sp\sixroman 1\sp,\ssp\sixroman 2}\,$, we get \œ$
y\ar 1=y\ar 2$ as follows. There are $B\sbi\iota$ and $\biit y\sbi\iota$ with
\œ$(\sp x\ssp,y\sbi\iota\ssp,B\sbi\iota\ssp,\biit y\sbi\iota)\in\Gamma\sp$. As
we have \œ$\biit y\sbi\iota\KN1\fvalue\emptyset=y\ar 0$ and \œ$
(\sp x\ssp,\biit y\sbi\iota\KN1\fvalue i\ssp,
           \biit y\sbi\iota\KN1\fvalue i\ssp\yplus)\in h$ for all \œ$
i\in\No\ssp$, since $h$ is a function, by induction we get \œ$
\bosy y\ar 1\KN1\fvalue i = \bosy y\ar 2\KN1\fvalue i$ for \œ$i\in\No\ssp$,
and hence \œ$\biit y\ar 1=\biit y\ar 2\,$. Since $\taurd E$ is a Hausdorff
topology,\linebreak and as we have $\biit y\ar 1\to y\sbi\iota$ in top $
\taurd E\,$, it follows that $y\ar 1=y\ar 2\,$.

We next prove that \œ$\dom g\ar 1=Q\ssp$. Trivially having \œ$\dom g\ar 1 \inc
Q\ssp$, arbitrarily given \œ$x\in x\ar 0+\tilde f\ssp'(\sp y\ar 0)\image B$
with \œ$B\in\Cal B\,$, it suffices to show that there are $y\ssp,\biit y$ with \œ$
(\sp x\ssp,y\ssp,B\ssp,\biit y\sp)\in\Gamma\sp$. To establish this, we
construct $\biit y$ by the following recursion:\linebreak fixing any \œ$z\ar 0
\in\Univ\snn\setminus\sn\vecs E\,$, we require that \œ$\biit y\fvalue\emptyset
=y\ar 0$ and \œ$\biit y\fvalue i\ssp\yplus =
h\fvalue(\sp x\ssp,\biit y\fvalue i\sp)$ in case\linebreak \œ$
(\sp x\ssp,\biit y\fvalue i\sp)\in\dom h\ssp$, otherwise putting \œ$
\biit y\fvalue i\ssp\yplus = z\ar 0\,$, for all \œ$i\in\No\ssp$. Letting $
(\sp l\sp)\ar A$ mean that \œ$\|\,\biit y\fvalue i\ssp\yplus\sn\yplus\snn -
     \biit y\fvalue i\ssp\yplus\ssp\|\sNor{X\sp B}\le\eps\,
\|\,\biit y\fvalue i\ssp\yplus\snn-
     \biit y\fvalue i\,\|\sNor{X\sp B}<\plusinfty$ holds for all $
i\in l\,\yplus\sp$, we first establish $\ssp\all{l\in\No}\,(\sp l\sp)\ar A$ by
induction as follows.

To get $(\emptyset)\ar A\,$, first note that $\biit y\fvalue 1\adot =
\biit y\fvalue\emptyset\ssp\yplus =
\ell\ssp\fvalue x-\ell\ssp\fvalue x\ar 0 + y\ar 0\,$. To get \œ$
(\sp x\ssp,\biit y\fvalue 1\adot)\in$ $\dom h\ssp$, by our arrangements and
Definitions \ref{Colo 0} it suffices for \œ$v\ar 0 =
\biit y\fvalue 1\adot-y\ar 0$ that $v\ar 0\in 2\,B\ssp$. This holds by \œ$
v\ar 0 = \ell\ssp\fvalue(\sp x-x\ar 0) \in
\ell\ssp\,[\,\sp(\sp\tilde f\ssp'(\sp y\ar 0))\image B\sp\,] = B\inc 2\,B\ssp
$, which also gives $\,\|\,v\ar 0\ssp\|\sNor{X\sp B}\le 1\ssp$. We hence have\vskip.3mm

$\mhyppy{6.7} \biit y\fvalue 2\adot =
\biit y\fvalue\emptyset\ssp\yplus\sn\yplus =
h\fvalue(\sp x\ssp,\biit y\fvalue 1\adot) =
\ell\ssp\fvalue x-\ell\sp\circ f\sp\fvalue(\sp\biit y\fvalue 1\adot)
  + \biit y\fvalue 1\adot\,$, \hfill whence\vskip.3mm\centerline{$
\biit y\fvalue 2\adot-\biit y\fvalue 1\adot =
\ell\ssp\fvalue x-\ell\sp\circ f\sp\fvalue(\sp\biit y\fvalue 1\adot) =
\ell\ssp\fvalue x\ar 0 + \biit y\fvalue 1\adot - y\ar 0 -
\ell\sp\circ f\sp\fvalue(\sp\biit y\fvalue 1\adot)$}\vskip.3mm

$\mhyppy{24.7} = v\ar 0 - \ell\ssp\fvalue(\sp
   f\sp\fvalue(\sp\biit y\fvalue 1\adot) - f\sp\fvalue y\ar 0)=
f\aar 1\KN1\fvalue y\ar 8 - f\aar 1\KN1\fvalue y\ar 0\,$,\vskip.3mm

\noin taking $y\ar 8=\biit y\fvalue 1\adot\,$.
Remark \ref{notes abou Colo 0}\ssp(c) now gives the assertion.

With $\sp l\in\No$ now assuming that $(\sp l\sp)\ar A$ holds, we prove $
(\sp l\ssp\yplus)\ar A$ as follows. Since we have $(\sp l\sp)\ar A\,$, for \œ$
v\ar 1 = \biit y\fvalue l\ssp\yplus\sn\yplus\snn - \biit y\fvalue l\ssp\yplus$
and \œ$v\ar 2 = \biit y\fvalue l\ssp\yplus\sn\yplus\sn\yplus\snn -
                 \biit y\fvalue l\ssp\yplus\sn\yplus$ we only have to
establish $\|\,v\ar 2\ssp\|\sNor{X\sp B}\le\eps\,\|\,v\ar 1\ssp\|\sNor{X\sp B}\,
$. By $(\sp l\sp)\ar A$ we have\vskip.5mm

$\mhyppy{13.2}
\|\,\biit y\fvalue l\ssp\yplus\sn\yplus\snn - y\ar 0\ssp\|\sNor{X\sp B} \le
\sum_{\,i\ssp\in\ssp l^{++\,}}\|\,\biit y\fvalue i\ssp\yplus\snn -
                       \biit y\fvalue i\,\|\sNor{X\sp B}$\nopagebreak\vskip.5mm\nopagebreak\noin
(e) $\mhyppy{38.7} \le \sum_{\,i\ssp\in\ssp l^{++\,}}\eps^{\sp\,i\,}
\|\,v\ar 0\ssp\|\sNor{X\sp B} \le \sum_{\,i\ssp\in\ssp l^{++\,}}\eps^{\sp\,i}
\le 2\,$,\vskip.5mm

\noin and consequently $(\sp x\ssp,\biit y\fvalue l\ssp\yplus\sn\yplus) \in
\dom h\,$, whence further\vskip.3mm

$\null\hfill   \biit y\fvalue l\ssp\yplus\sn\yplus\sn\yplus =
\ell\ssp\fvalue x-\ell\sp\circ f\sp\fvalue(\sp
   \biit y\fvalue l\ssp\yplus\sn\yplus) +
      \biit y\fvalue l\ssp\yplus\sn\yplus\sp$. \hyppy{23mm} Also having\vskip.3mm

$\null\hfill   \biit y\fvalue l\ssp\yplus\sn\yplus =
   \ell\ssp\fvalue x-\ell\sp\circ f\sp\fvalue(\sp
     \biit y\fvalue l\ssp\yplus) + \biit y\fvalue l\ssp\yplus\sp$, \hyppy{35.3mm}
we get\vskip.3mm

\centerline{$v\ar 2 = v\ar 1 - \ell\ssp\fvalue(\sp f\sp\fvalue(\sp
  \biit y\fvalue l\ssp\yplus\sn\yplus) -
    f\sp\fvalue(\sp\biit y\fvalue l\ssp\yplus))=
    f\aar 1\KN1\fvalue(\sp\biit y\fvalue l\ssp\yplus\sn\yplus) -
    f\aar 1\KN1\fvalue(\sp\biit y\fvalue l\ssp\yplus)\,$,}\vskip.3mm

\noin whence again Remark \ref{notes abou Colo 0}\ssp(c) gives the assertion.

Now having obtained \œ$\ssp\all{l\in\No}\,(\sp l\sp)\ar A\,$, we know that $
\biit y$ is a Cauchy sequence in $X\sbi B\,$, which is Banachable, hence
complete by the assumption that $E$ is locally complete. Hence, there is $y$
with \œ$\biit y\to y$ in top $\taurd(X\sbi B\spp)\,$, and consequently also $
\biit y\to y$ in top $\taurd E\,$. We have now concluded the proof of $
\dom g\ar 1=Q\ssp$.

From the assumption that $\Cal B$ is a bornological generator for $E$ it
follows that $\bigcup\ssp\Cal B$ is an absolutely convex bornivore in $E\,$,
hence also $\sp\tilde f\ssp'(\sp y\ar 0)\image\bigcup\ssp\Cal B$ since $\sp
\tilde f\ssp'(\sp y\ar 0)$ is assumed to be a linear homeomorphism \œ$E\to E\,
$. Since $E$ is assumed to be bornological, we have \œ$\ssp
\tilde f\ssp'(\sp y\ar 0)\image\bigcup\ssp\Cal B\in\ymp E\,$, whence it
follows existence of some \œ$\ssp U\in\taurd E$ with \œ$x\ar 0\in U \inc
Q\ssp$. Then putting \œ$g=g\ar 1\ssp|\,U\spp$, we have $g$ a function
with $\sp x\ar 0\in\dom g=U\in\taurd E\sp\cap\Pows Q\ssp$, recalling that
$\Pows Q=\{\,S:S\inc Q\,\}\,$.

We now proceed to prove \œ$\ssp g\inc f\inve\sp$. For this considering
arbitrary \œ$(\sp x\ssp,y\sp)\in g\,$, we have \œ$(\sp x\ssp,y\sp) \in
g\ar 1\ssp$, whence there are $B\ssp,\biit y$ with \œ$
(\sp x\ssp,y\ssp,B\ssp,\biit y\sp)\in\Gamma\sp$. A slight rearrangement of the
arguments used to establish \œ$\dom g\ar 1=Q$ shows that (e) holds. This
gives \œ$\biit y\in(\sp y\ar 0+2\,B\sp)\potNo$ whence by closedness of $B$ we
get \œ$y\in y\ar 0+2\,B\inc\dom f\sp$.\linebreak Since for \œ$P =
(\sp\taurd(X\sbi B\spp)\ssp,\taurd E\ssp)$ we have \œ$(\spp P\spp,f\,|\,
(\sp y\ar 0+2\,B\ssp))$ a continuous map, also $(\spp P\spp,\sp\ell\sp\circ
f\sp\,|\,(\sp y\ar 0+2\,B\ssp))$ is such. For \œ$i\in\No$ having \œ$
\biit y\fvalue i\ssp\yplus =
\ell\ssp\fvalue x-\ell\sp\circ f\sp\fvalue(\sp\biit y\fvalue i\sp)
  + \biit y\fvalue i\ssp$, we get $y = 
\ell\ssp\fvalue x-\ell\sp\circ f\sp\fvalue y + y\,$, \,consequently $
  x=f\sp\fvalue y\ssp$, and hence $(\sp x\ssp,y\sp)\in f\inve\sp$.

Fixing \œ$\ssp B\in\Cal B$ and \œ$\,x\ar 1\ssp,x\ar 2 \in(\sp x\ar 0 +
\sp\ell\ssp\inve\image\snn B\ssp)\spp\cap U\spp$, for \œ$y\sbi\iota =
g\fvalue x\sbi\iota$ and \œ$u = x\ar 1\snn-x\ar 2$ and
\œ$v=y\ar 1\snn-y\ar 2$ it remains to establish \œ$\,
\|\ssp v\ssp\|\sNor{X\sp B} \le
(\sp 1-\eps\sp)^{\sp\mminus 1}\sp\|\,\sp\ell\ssp\fvalue u\,\|\sNor{X\sp B}\,$.
We have \œ$y\sbi\iota = \ell\ssp\fvalue x\sbi\iota -
\ell\sp\circ f\sp\fvalue y\sbi\iota+ y\sbi\iota\,$, \,whence we get\vskip.3mm

$\mhyppy8  v-\ell\ssp\fvalue u =
v-\ell\ssp\fvalue(\sp f\sp\fvalue y\ar 1\snn - f\sp\fvalue y\ar 2) =
f\aar 1\KN1\fvalue y\ar 1\snn - f\aar 1\KN1\fvalue y\ar 2\,$. \hfill
Noting that by the \vskip.3mm\noin preceding paragraph we have $y\ar 1\ssp,\sp
y\ar 2\in y\ar 0+2\,B\ssp$, by Remark \ref{notes abou Colo 0}\ssp(c) we get\vskip.5mm$\mhyppy{6.5}
\|\,v-\ell\ssp\fvalue u\,\|\sNor{X\sp B} =
\|\, f\aar 1\KN1\fvalue y\ar 1\snn -
     f\aar 1\KN1\fvalue y\ar 2\ssp\|\sNor{X\sp B} \le
\eps\,\|\ssp v\ssp\|\sNor{X\sp B}\,$\vskip.5mm

\noin and further \ \ $\|\ssp v\ssp\|\sNor{X\sp B} =
\|\,v-\ell\ssp\fvalue u+\ell\ssp\fvalue u\,\|\sNor{X\sp B}$\vskip.3mm

$\mhyppy{28.3} \le \|\,v-\ell\ssp\fvalue u\,\|\sNor{X\sp B} +
  \|\,\sp\ell\ssp\fvalue u\,\|\sNor{X\sp B} \le
\eps\,\|\ssp v\ssp\|\sNor{X\sp B} +
     \|\,\sp\ell\ssp\fvalue u\,\|\sNor{X\sp B}\,$,\vskip.5mm

\noin whence finally \ \ $
\|\ssp v\ssp\|\sNor{X\sp B} \le (\sp 1-\eps\sp)^{\sp\mminus 1}\sp
  \|\,\sp\ell\ssp\fvalue u\,\|\sNor{X\sp B}\,$.                  \end{proof}

\begin{theorem}\label{inve func}

If \œ$\,k\not=\emptyset$ and \œ$\,\tilde f=(E\ssp,E\ssp,f\sp)\in\Lip_FKt^k$
and $\,\tilde f$ is Colombeau$\ar{\,0\,}$--\,tame at $y\ar 0\,,$ there is $\ssp
U$ with $\ssp y\ar 0\in U\in\tauMac E\sp$ and $\,(E\ssp,E\ssp,
(\sp f\,|\,U\ssp)\inve\spp)\in\Lip_FKt^k\,$.                   \end{theorem}

\begin{proof} Assuming the premise, let \œ$x\ar 0=f\sp\fvalue y\ar 0$ and \œ$
X = \sigrd E\,$. Now, there is $\Cal B\ar 0$ such that \ssp Colo$\ar{\,0}\ssp
(\ssp\tilde f\,|\subtext{map\,}\vecs E\,,\sp y\ar 0\,,\sp\frac 12\,,
\Cal B\ar 0)\ssp$ holds. Putting \œ$Q\ar 0=y\ar 0+2\,\bigcup\ssp\Cal B\ar 0\,
$, we show indirectly that $f\,|\,Q\ar 0$ is injective. Indeed, if this does
not hold, there are distinct \œ$y\ar 1\ssp,y\ar 2\in\dom\sn f\sp\cap Q\ar 0$
with \œ$f\sp\fvalue y\ar 1 = f\sp\fvalue y\ar 2\,$, and we find some \œ$B \in
\Cal B\ar 0$ with \œ$\sp y\ar 1\ssp,\sp y\ar 2\in y\ar 0+2\,B\ssp$. Letting \œ$
v = y\ar 1\snn-y\ar 2\,$, by Remark \ref{notes abou Colo 0}\ssp(c) we then get
\œ$\|\ssp v\ssp\|\sNor{X\sp B}=$\Biggerlineskip1 \œ$
\|\,y\ar 1\snn - y\ar 2\ssp\|\sNor{X\sp B} =
 \|\,f\aar 1\KN1\fvalue y\ar 1\snn -
     f\aar 1\KN1\fvalue y\ar 2\ssp\|\sNor{X\sp B} \le
\frac 12\,\|\,y\ar 1\snn - y\ar 2\ssp\|\sNor{X\sp B} =
\frac 12\,\|\ssp v\ssp\|\sNor{X\sp B}$ whence finally\Biggerlineskip1 $\,
0 < \|\ssp v\ssp\|\sNor{X\sp B} =
2\,\|\ssp v\ssp\|\sNor{X\sp B} - \|\ssp v\ssp\|\sNor{X\sp B} \le
\|\ssp v\ssp\|\sNor{X\sp B} -\|\ssp v\ssp\|\sNor{X\sp B} = 0\,$, \sp a
contradiction. \vskip.3mm

Next, since $\ssp\bigcup\ssp\Cal B\ar 0$ is an absolutely convex bornivore in
$E\,$, there is \œ$\ssp V\in\taurd E$ with \œ$y\ar 0\in V\inc Q\ar 0\,$,
whence there further is some $\Cal B$ with \ssp Colo$\ar{\,0}\ssp(\ssp
\tilde f\,|\subtext{map\,}V\spp,\sp y\ar 0\,,\sp\frac12\,,\Cal B\ssp)\,$.\Biggerlineskip1
Letting $g$ be as given by Lemma \ref{on impl fus} above, and taking \œ$\sp
U = \rng g\,$, since $\sp f\sp\,|\,\sp V\sp$ is injective and \œ$\sp g \inc
(\sp f\sp\,|\,\sp V\ssp)\inve\sp$, and as \œ$\dom g\in\taurd E\inc\tauMac E\,
$, in view of continuity of $(\sp\tauMac E\ssp,\tauMac E\ssp,
f\sp\,|\,\sp V\ssp)$ it follows that \œ$\sp U =
(\sp f\sp\,|\,\sp V\ssp)\inve\sp[\,\dom g\sp\,]\in\tauMac E\,$. Trivially
having $\ssp y\ar 0\in U\spp$, for $\sp\tilde g = (E\ssp,E\ssp,g\sp)$ it
remains to establish $\sp \tilde g\in\Lip_FKt^k\,$.

For this using \cite[Theorem 4.8.4, p.\ 152]{FK}\ssp, in view of
Remark \ref{notes abou Colo 0}\ssp(\sp b) it suffices that \œ$\tilde g \in
\Lip_FKt^{0.}\,$. That is, for arbitrarily given \œ$c \in
\Lip_FKt^{0.}\KN2\image\{\ssp(\biit R\,,E\ssp)\ssp\}$ and \œ$\gamma=g\circ c\,
$, we should have \œ$\gamma \in
\Lip_FKt^{0.}\KN2\image\{\ssp(\biit R\,,E\ssp)\ssp\}\,$. To get this, for
arbitrarily given \œ$t\ar 0\in\dom\gamma$ it suffices to show existence of \œ$
\delta\in\Rep$ and \œ$B\in\Cal B$ such that for \œ$J =
{\,]\sp}\,t\ar 0-\delta\ssp,t\ar 0+\delta\,{\sp[\,\sp}$ we have \œ$\sp
\gamma\sp\,|\,J\in\Lip_FKt^{0.}\KN2\image\{\ssp(\biit R\,,X\sbi B\spp)\ssp\}\,
$. To get this, we put \œ$\Cal B\ar 1 = \sp\ell\ssp\inve\images\Cal B\ssp$,
and first note that \œ$\sp c\fvalue t\ar 0-y\ar 0\in\dom g-y\ar 0\in\taurd E\,
$. Hence, there is a real \œ$r>1$ with \œ$r\,(\sp c\fvalue t\ar 0-y\ar 0)$ \œ$
\in\dom g-y\ar 0\inc\bigcup\ssp\Cal B\ar 1\ssp$, whence further there is \œ$
B\ar 2\in\Cal B\ar 1$ with \œ$r\,(\sp c\fvalue t\ar 0-y\ar 0)\in B\ar 2\,$. By
\cite[Corollary 1.8, p.\ 13]{KM} we then find \œ$\delta\in\Rep$ and \œ$
B\ar 1\in\Cal B\ar 1$ with \œ$B\ar 2\inc B\ar 1$ and \œ$\sp c\sp\,|\,J \in
\Lip_FKt^{0.}\KN2\image\{\ssp(\biit R\,,X\sbi{B_1})\ssp\}\,$. Since now \œ$
c\fvalue t\ar 0-y\ar 0 \in
\roman{Int}\sbi{\,\eightmath\tauu{_{\!}}_{rd}\sp(X_{B_1})\,}B\ar 1\ssp$, we
may take $\delta$\Biggerlineskip1 smaller so that also \œ$c\ssp\image\sn J\inc
y\ar 0 + B\ar 1\ssp$. For \œ$B=\sp\ell\ssp\image\snn B\ar 1\ssp$, the \q{in
addition} part of Lemma \ref{on impl fus} now gives $\sp\gamma\sp\,|\,J \in
\Lip_FKt^{0.}\KN2\image\{\ssp(\biit R\,,X\sbi B\spp)\ssp\}\,$.   \end{proof}

\begin{corollary}\label{Frec inve fu}

If \œ$\sp\tilde f=(E\ssp,E\ssp,f\sp)\in\C_c^\infty(R)$ with $E$ a
cm-convenient space$\ssp,$ and if also $\sp\tilde f$ is Colombeau$\ar{\,0\,}
$--\,tame at $x\ar 0\,,$ there is $U$ with\par\centerline{$
x\ar 0\in U\in\taurd E\sp$ and $\,(E\ssp,E\ssp,
(\sp f\,|\,U\ssp)\inve\spp)\in\C_c^\infty(R)\,$.}            \end{corollary}

\begin{proof} Since \œ$\sp\tauMac E=\taurd E\ssp$ holds for cm-convenient $E\,
$, the assertion immediately follows from
Proposition \ref{rela wtw Lip^infty and C_c^infty} and
Theorem \ref{inve func} above.                                   \end{proof}

To indicate the basic idea for proving Colombeau$\ar{\,0\,}$--\,tameness for
maps of the type \œ$\,x\mapsto\varphi\circ[\,\roman{id\,},x\,]\,$, and also to
show that Colombeau$\ar{\,0\,}$--\,tameness despite of its strength is not too
restrictive, we establish the following

\begin{proposition}\label{appl exa}

Let \œ$\ssp I=[\,0\,,1\,]$ and \œ$\,E=\Cinfty(I\sp)\,,$ and also let \œ$\ssp
\varphi:I\timesn\Re\to\Re$ be smooth with \œ$\ssp 0 \not\in
\rng\partial\ar 2\ssp\varphi\,$. For \œ$\ssp f=\seq{\,\varphi\circ
[\,\roman{id\,},x\,]:x\in\vecs E\sp\,}$ and for \œ$\ssp\tilde f =
(E\ssp,E\ssp,f\sp)\,,$ then $\tilde f$ is Colombeau$\ar{\,0\,}$--\,tame at
every $\,x\in\vecs E\,$.                                   \end{proposition}

\begin{proof} Assuming that \œ$x\in V\in\taurd E\,$, there is \œ$l\ar 0\in\No$
with the property that $x+V\aar 0\inc V$ for the set $\ssp V\aar 0 =
\vecs E\sp\cap\{\,z:\all{i\in l\ar 0\KN1\yplus\sp,s\in I}\,
l\ar 0\,|\,z^{(i)}\fvalue s\,|\le 2\,\}\,$. For\vskip1mm

$\mhyppy{5}\chi=I\timesn\Re\times\snn\Re\ssp\cap\big\{\,
(\sp s\ssp,\eta\ssp,t\sp):$

$\mhyppy{21} t = 4\,
(\sp\partial\ar 2\ssp\varphi\fvalue(\sp s\ssp,x\fvalue s\sp)
)^{\sp\mminus 1}\spp\int_{\,0}^{\,1}
\partial\ar 2^{\,2.}\ssp\varphi\fvalue(\sp s\ssp,x\fvalue s+2\,
s\ar 1\ssp\eta\sp)\,\ssp\roman d\,s\ar 1\ssp\big\}\,$,\vskip1mm

\noin we observe that \œ$\,2\,(\sp\tilde f\ssp'(x))\inve\snn\circ
(\sp\tilde f\ssp'(\sp x+2\,u\sp)-\tilde f\ssp'(x))\fvalue v =
\chi\circ[\,\sp\roman{id\,},u\,]\cdotn u\cdotn v\,$ for any\linebreak \œ$
u\ssp,v\in\vecs E\,$. To prove that $\tilde f$ is Colombeau$\ar{\,0\,}
$--\,tame at $x\ssp$, it hence suffices to establish a bornological generator
$\Cal B$ for $E$ such that we have $\ssp x\in\vecs(\sp\sigrd E\sbi B)$ and \vskip.5mm

\noin\ ($*$) $\mhyppy{15}
\chi\circ[\,\sp\roman{id\,},u\,]\cdotn u\cdotn v\in B\inc
\frac 12\,V\aar 0$ \ whenever \ $u\ssp,v\in B\in\Cal B\ssp$.\vskip.5mm

To get this, we make the following preparations and observations. Write\vskip.3mm

\centerline{$\roman B\,\biit m = \vecs E\sp\cap\{\,x:
\all{i\in\No\ssp, s\in I}\,|\,x^{(i)}\fvalue s\,|\le\bosy m\fvalue i\,\}\,$,}\vskip.3mm

\noin where we generally require \œ$\biit m\in(\Rep)\potNo$ to be
nondecreasing. For \œ$i\in\No\ssp$, let \œ$\ssp\roman P\ssp i =
\Nopot{\times 2.}\cap\{\ssp(\sp i\ar 1\sp,i\ar 2) :
i\ar 1\snn + i\ar 2\in i\ssp\yplus\ssp\}\,$, and for \œ$\chi\ar 1 \in
\vecs\Cinfty\spp(\spp I\timesn\Re\sp)$ and \œ$u\in\vecs E$\linebreak define
the \q{jet} functions\vskip.5mm

\centerline{$\ssp\roman J\ar 2^{\,i}\sp\chi\ar 1\sn:I\timesn\Re\owns
\zeta\mapsto \{\ssp(\sp i\ar 1\sp,i\ar 2\ssp,
\partial_{\sp\sixroman 1}^{\,\sp i_1}\sp
\partial_{\sp\sixroman 2}^{\,\sp i_2}\sp\chi\fvalue\zeta\sp):
(\sp i\ar 1\sp,i\ar 2)\in\roman P\ssp i\,\}
\in\Re\,^{\sp\roman P\sp i}$}\vskip.3mm

\noin and  $\mhyppy{9.5} \roman J\ar 1^{\,i}\sp u:I\owns
s \mapsto \seq{\,u^{(l)}\fvalue s:l\in i\ssp\yplus\ssp}
\in\Re\,^{i\ssp+\sp 1.}\sp$.\vskip.5mm

If we have a polymial function \œ$\ssp
p : \Re\,^{\sp\roman P\sp i}\timesn\Re\,^i\timesn\Re\,^i\to\Re$ with the
properties\Biggerlineskip1 that for fixed $\bosy\xi$ the map \œ$
(\sp\bosy\eta\ssp,\bosy\zeta\sp) \mapsto
  p\fvalue(\sp\bosy\eta\ssp,\bosy\zeta\ssp,\bosy\xi\sp)$ is bilinear \œ$\ssp
\bold R\expnota^\roman P\sp i]_{vs}\vstimes(\ssp\bold R\expnota^i]_{vs}\snn)
\to\bold R\sp\,$, and for all $\chi\ar 1 \in
\vecs\Cinfty\spp(\spp I\timesn\Re\sp)$ and $u\ssp,v\in\vecs E$ we have\vskip.6mm

$\mhyppy{5}
(\sp\chi\ar 1\snn\circ[\,\sp\roman{id\,},u\,]\cdotn u\cdotn v\sp)^{(i)} =
\partial\ar 2\ssp\chi\ar 1\snn\circ[\,\sp\roman{id\,},u\,]\cdotn u^{(i)}
\cdotn u\cdotn v + \chi\ar 1\snn\circ[\,\sp\roman{id\,},u\,]\cdotn u^{(i)}
\cdotn v$\nopagebreak\par\nopagebreak\noin
(s) $\mhyppy{19} + \sp\chi\ar 1\snn\circ[\,\sp\roman{id\,},u\,]\cdotn u\cdotn
v^{(i)} + p\circ[\,\sp
\roman J\ar 2^{\,i\,}\chi\ar 1\snn\circ[\,\sp\roman{id\,},u\,]\,,\sp
\roman J\ar 1^{\,i\ssp-\sp 1.\,}v\ssp,\sp
\roman J\ar 1^{\,i\ssp-\sp 1.\,}u\sp\,]\,$,\vskip.6mm

\noin it follows that \ \ \ $(\sp\chi\ar 1\snn\circ
         [\,\sp\roman{id\,},u\,]\cdotn u\cdotn v\sp)^{(i\ssp+\sp 1.)} = $

$\mhyppy6   \partial\ar 1\ssp\partial\ar 2\ssp\chi\ar 1\snn\circ
[\,\sp\roman{id\,},u\,]\cdotn u^{(i)}\cdotn u\cdotn v +
\partial\ar 2^{\,\sp 2.}\ssp\chi\ar 1\snn\circ[\,\sp\roman{id\,},u\,]
\cdotn u\ssp'\snn\cdotn u^{(i)}\cdotn u\cdotn v $

$\mhyppy8   + \partial\ar 2\ssp\chi\ar 1\snn\circ[\,\sp\roman{id\,},u\,]\cdotn
      u^{(i\ssp+\sp 1.)}\cdotn u\cdotn v + \partial\ar 2\ssp\chi\ar 1\snn\circ
         [\,\sp\roman{id\,},u\,]\cdotn u^{(i)}\cdotn u\ssp'\snn\cdotn v $

$\mhyppy8   + \partial\ar 2\ssp\chi\ar 1\snn\circ[\,\sp\roman{id\,},u\,]\cdotn
           u^{(i)}\cdotn u\cdotn v\ssp' + \partial\ar 1\ssp\chi\ar 1\snn\circ
            [\,\sp\roman{id\,},u\,]\cdotn u^{(i)}\cdotn v$ 

$\mhyppy8   + \partial\ar 2\ssp\chi\ar 1\snn\circ[\,\sp\roman{id\,},u\,]\cdotn
               u\ssp'\snn\cdotn u^{(i)}\cdotn v + \chi\ar 1\snn\circ
                [\,\sp\roman{id\,},u\,]\cdotn u^{(i\ssp+\sp 1.)}\cdotn v$ 

$\mhyppy8   + \chi\ar 1\snn\circ[\,\sp\roman{id\,},u\,]\cdotn u^{(i)}\cdotn
               v\ssp' + \partial\ar 1\ssp\chi\ar 1\snn\circ
                 [\,\sp\roman{id\,},u\,]\cdotn u\cdotn v^{(i)}$ 

$\mhyppy8   + \partial\ar 2\ssp\chi\ar 1\snn\circ[\,\sp\roman{id\,},u\,]\cdotn
               u\ssp'\snn\cdotn u\cdotn v^{(i)} + \chi\ar 1\snn\circ
                 [\,\sp\roman{id\,},u\,]\cdotn u\ssp'\cdotn v^{(i)}$ 

$\mhyppy8   + \chi\ar 1\snn\circ[\,\sp\roman{id\,},u\,]\cdotn u\cdotn
              v^{(i\ssp+\sp 1.)} +
     p\circ[\,\sp\roman J\ar 2^{\,i\,}\partial\ar 1\ssp\chi\ar 1\snn\circ
[\,\sp\roman{id\,},u\,]\,,\sp\roman J\ar 1^{\,i\ssp-\sp 1.\,}v\spp,\sp
                             \roman J\ar 1^{\,i\ssp-\sp 1.\,}u\sp\,]$

$\mhyppy8 +
p\circ[\,\sp\roman J\ar 2^{\,i\,}\partial\ar 2\ssp\chi\ar 1\snn\circ
[\,\sp\roman{id\,},u\,]\,,\sp\roman J\ar 1^{\,i\ssp-\sp 1.\,}v\spp,\sp
                       \roman J\ar 1^{\,i\ssp-\sp 1.\,}u\sp\,]\cdotn u\ssp'$

$\mhyppy8 +
p\circ[\,\sp\roman J\ar 2^{\,i\,}\chi\ar 1\snn\circ[\,\sp\roman{id\,},u\,]\,,
\sp\roman J\ar 1^{\,i\ssp-\sp 1.\,}v\ssp'\spp,\sp
   \roman J\ar 1^{\,i\ssp-\sp 1.\,}u\sp\,]$

$\mhyppy8 +
\partial\ar 3\ssp p\circ[\,\sp\roman J\ar 2^{\,i\,}\chi\ar 1\snn\circ
[\,\sp\roman{id\,},u\,]\,,\sp\roman J\ar 1^{\,i\ssp-\sp 1.\,}v\ssp,\sp
\roman J\ar 1^{\,i\ssp-\sp 1.\,}u\sp\,]\,.\,(\sp
\roman J\ar 1^{\,i\ssp-\sp 1.\,}u\ssp'\spp)\,$.\vskip1mm

Omitting the details, from the preceding one sees that by a suitable recursion
one can construct \œ$\,\bold P\in\Univ\,^{I\!\!N}$ such that for \œ$
(\sp i\ssp,p\ssp)\in\bold P$ we have \œ$p=\bold P\sp\fvalue i$ a polynomial as
above such that (s) holds for the appropriate $\ssp\chi\ar 1\ssp,\sp u\ssp,\sp
v\ssp$. Letting $\,\roman R\,i\,s$ denote\vskip.3mm

$\mhyppy{5} \sup\ssp\{\,|\,\sp\bold P\sp\fvalue i\ssp\yplus\fvalue
(\sp\bosy\xi\,,\bosy\eta\,,\bosy\zeta\ssp)\,|:i:\bosy\xi\in
\roman J\ar 2^{\,i\ssp+\sp 1.\,}\chi\,[\,I^{\sp\times 2.}\ssp]$

\hyppy{48mm} and $\sup\ssp\{\,|\ssp r\sp|:
r\in\rng\snn(\sp\bosy\eta\sp\cup\sp\bosy\xi\ssp)\ssp\}\le s\,\}\,$,

\noin and constructing $\sp\bosy\rho\in\Univ\potNo$ by the recursion\vskip.6mm

$\mhyppy{17.6}\bosy\rho\fvalue\emptyset \ \, = \seq{\,\max\,\{\ssp s\ssp,
\roman R\,\emptyset\,s\ssp\}:s\in\Repp\ssp}\mhyppy{16}$ and

\centerline{$\bosy\rho\fvalue i\ssp\yplus = \seq{\,\max\,\{\,\sp
\bosy\rho\fvalue i\fvalue s\,,\sp\roman R\,i\ssp\yplus\ssp s\,\}:s \in
\Repp\ssp}$ \ for \ $i\in\No\ssp$,}\vskip.6mm

\noin and putting \œ$\sp\rho=\bosy\rho\sp\LHB{.9}{^{^{\,\wedge}}}\sp$, then $\sp
\rho$ is a function \œ$\No\timesn\Repp\to\Repp\,$, nondecreasing separately in
both arguments, with \œ$\rho\fvalue(\sp i\ssp,0\sp)=0$ and \œ$s \le
\rho\fvalue(\sp i\ssp,s\sp)$ for \œ$i\in\No$ and \œ$s\in\Repp\,$, and $
\rho\,(\sp i\ssp,\cdot\sp)$ continuous at $0\,$, and also such that with

$\mhyppy{13}\smb B\ar 0=\sup\ssp\{\,1+|\ssp r\sp|:r \in
(\sp\chi\cup\partial\ar 2\ssp\chi\sp)\,[\,I^{\sp\times 2.}\ssp]\sp\,\}
\mhyppy{14}$ we have\vskip.5mm\centerline{$
|\,(\sp\chi\circ[\,\sp\roman{id\,},u\,]\cdotn u\cdotn v\sp)^{(i\ssp+\sp 1.)}
\fvalue s\,| \le \smb B\ar 0\,(\sp\smb M\ar 0+2\sp)\,\smb M\ar 0\,
(\sp\bosy m\fvalue i\ssp\yplus) + \rho\fvalue(\sp i\ssp,\smb M\sp)$}\vskip.5mm

\noin for \œ$i\in\No$ and \œ$s\in I$ and \œ$u\ssp,v\in\roman B\,\bosy m$ and \œ$
\bosy m\in(\Rep)\potNo$ such that \œ$\bosy m\fvalue\emptyset\le\smb M\ar 0
\le 1$ and $\ssp\sup\,(\sp\bosy m\image i\ssp\yplus)\le\smb M<\plusinfty\,$.

Putting \ \ \ $\theta\ar 0 = \Rep\ssn\times\Repp\sn\times\No\times\Repp\cap
                         \{\ssp(\sp r\sp,s\ssp,i\ssp,t\sp):$

$\mhyppy{30} \rho\fvalue(\sp i\ssp,s\sp) = t\,(\sp 1-\smb B\ar 0\,r\,(\sp 2 +
r\sp))$ and $\smb B\ar 0\,r\,(\sp 2+r\sp)<1\,\}\,$,\vskip.5mm

\noin we note that $\theta\ar 0\sp(\sp r\sp,\ssp\cdot\,,i\sp)$ is continuous
at $0$ with \œ$\theta\ar 0\KN1\fvalue(\sp r\sp,0\ssp,i\sp)=0$ whenever \œ$
(\sp r\sp,0\ssp,i\sp)\in$ $\dom\theta\ar 0\,$. Using this, by a suitable
finite induction, one first establishes existence of \œ$\bosy n \in
(\Rep)\,^{l_0\sp+\sp 1.}$ with \œ$\{\ssp(\sp\bosy n\fvalue\emptyset\,,
\bosy n\fvalue i\ssp,i\ssp,\bosy n\fvalue i\ssp\yplus):i\in l\ar 0\ssp\} \inc
\theta\ar 0$ and \œ$\bosy n\fvalue\emptyset \le
\frac 13\,\smb B_{\sp\sixroman 0}^{\ssp\mminus 1}$ and $\ssp
l\ar 0\,(\sp\bosy n\fvalue l\ar 0)\le 1\ssp$, and then fixes one such $\sp
\bosy n\,$. With\vskip.5mm

$\mhyppy{2.65} \biit x\ar 0 = \seq{\,\sp\sup\,\{\,1+|\,x^{(l)}\fvalue s\,|:
s\in I$ and $\ssp l\in i\ssp\yplus\ssp\}:i\in\No\ssp}$ \hfill and\KP6\vskip.3mm

$\mhyppy{2.65} \biit x\ar 1 =
\seq{\,(\sp\biit n\fvalue\emptyset\sp)^{\sp\mminus 1}\sp
(\sp\biit x\ar 0\KN1\fvalue l\ar 0)^{\sp\mminus 1}\sp
(\sp\biit x\ar 0\KN1\fvalue i\ssp):i\in\No\ssp}$ \hfill and\KP6\vskip.3mm

$\mhyppy{4.75} \theta = \{\ssp(\sp r\sp,s\ssp,i\ssp,t\sp):\exi{t\ar 1}\,
(\sp r\sp,s\ssp,i\ssp,t\ar 1)\in\theta\ar 0$ and $\sp t = \max\,\{\ssp
t\ar 1\ssp,\sp\biit x\ar 1\KN1\fvalue i\ssp\yplus\ssp\}\sp\}$ \hfill and\KP6\vskip.3mm

\centerline{$\Cal M=(\Rep)\potNo\cap\{\,\bosy m:
\bosy n\inc\bosy m\text{ and }\ssp\all{i\in\No}\,
\theta\sp\fvalue(\sp\bosy m\fvalue\emptyset\,,\bosy m\fvalue i\ssp,i\ssp)
\le\bosy m\fvalue i\ssp\yplus\ssp\}\,$,}\vskip.3mm

\noin we now take $\Cal B=\{\,\sp\roman B\,\biit m:\biit m\in\Cal M\sp\,\}\,$.

Note that by \œ$\bosy m\fvalue i\le\rho\fvalue(\sp i\ssp,\bosy m\fvalue i\sp)
\le \theta\sp\fvalue(\sp\bosy m\fvalue\emptyset\,,\bosy m\fvalue i\ssp,i\sp)
\le\bosy m\fvalue i\ssp\yplus$ every \œ$\bosy m\in\Cal M$ is nondecreasing.
By our construction, it is straightforward to verify that ($*$) above holds.
Since every $B\in\Cal B$ is absolutely convex, to have $\Cal B$ a bornological
generator for $E\,$, one should verify \vskip.5mm

\noin\ (\sp i$\ar 1$) \ \ $
 \all{\bosy m\ar 1\ssp,\bosy m\ar 2\in\Cal M}\,\exi{\bosy m\in\Cal M}\,
  \all{i\in\No}\,\max\,\{\,\bosy m\ar 1\KN1\fvalue i\ssp,
                     \bosy m\ar 2\KN1\fvalue i\,\}\le\bosy m\fvalue i\,$,\vskip.3mm

\noin\ (\sp i$\ar 2$) \ \ $
 \all{\bosy b\in(\Rep)\potNo}\,\exi{\eps\in\Rep\sp,\bosy m\in\Cal M}\,
  \all{i\in\No}\,\eps\,(\sp\bosy b\fvalue i\sp)\le\bosy m\fvalue i\,$.\vskip.5mm

\noin To get (\sp i$\ar 1$)\ssp, for $i\in\No\sn\setminus\snn l\ar 0$ one
applies the recursion\vskip.3mm\centerline{$
\bosy m\fvalue i\ssp\yplus = \max\,\{\,
\bosy m\ar 1\KN1\fvalue i\ssp\yplus\sp,\sp
\bosy m\ar 2\KN1\fvalue i\ssp\yplus\sp,\sp
\theta\sp\fvalue(\sp\bosy n\fvalue\emptyset\,,\bosy m\fvalue i\ssp,i\ssp)\,\}
$ \ with \ $\bosy m\fvalue l\ar 0=\bosy n\fvalue l\ar 0\,$.}\vskip.3mm

\noin For (\sp i$\ar 2$) one first chooses \œ$\eps\in\Rep$ so that \œ$\ssp
l\ar 0\KN1\yplus\inc\{\,i:\eps\,(\sp\bosy b\fvalue i\sp) \le
\bosy n\fvalue i\,\}\,$, and then\linebreak applies the recursion \ $
\bosy m\fvalue i\ssp\yplus =
\max\,\{\,\eps\,(\sp\bosy b\fvalue i\ssp\yplus)\ssp,
\theta\sp\fvalue(\sp\bosy n\fvalue\emptyset\,,\bosy m\fvalue i\ssp,i\ssp)
\,\}\,$.                                                         \end{proof}

To give some perspective, we conclude with the following

\begin{remark}

Original formulations of Colombeau's \q{tameness} conditions are given in
\cite{Colo} and reproduced in \cite[Section \erm{XIII\spp}.4]{HN1}\ssp. There
these conditions concern Silva differentiable maps between convex bornological
vector spaces. In \cite{Les}\ssp, the conditions are adapted for maps between
bornological locally convex spaces. Definition 2.1 in \cite[p.\ 428]{Les}
introduces certain order $k$ differentiability classes which by
\cite[Theorem 2.8.1\sp(2)\ssp, pp.\ 102, 105]{Ke} and by a suitable adaptation
of the idea in the proof of \cite[Theorem 5.20, pp.\ 62, 27]{KM} are precisely
the classes $\C_c^k(R)\,$. It seems that in \cite{Les} one has not taken into
account carefully enough the fact that a Silva $\sp C^{\ssp 1.}$ map need not
be continuous with respect to the locally convex topologies when the domain
space is not cm-convenient.                                     \end{remark}


\end{document}